\documentclass[11pt,oneside,fleqn]{article}

\usepackage[ansinew]{inputenc}
\usepackage[mathscr]{eucal}
\usepackage{amsmath,amssymb,amsthm}
\usepackage{graphicx}
\usepackage{cite}
\usepackage{comment}
\usepackage{hyperref}
\usepackage{subcaption}
\allowdisplaybreaks

\hypersetup{colorlinks, linkcolor=blue, citecolor=blue, urlcolor=blue}
\makeatletter
\long\def\@makecaption#1#2{%
  \vskip\abovecaptionskip\footnotesize
  \sbox\@tempboxa{#1. #2}%
  \ifdim \wd\@tempboxa >\hsize
    #1. #2\par
  \else
    \global \@minipagefalse
    \hb@xt@\hsize{\hfil\box\@tempboxa\hfil}%
  \fi
  \vskip\belowcaptionskip}
\makeatother

\newcommand{\todo}[1][\null]{\ensuremath{\clubsuit}}

\newcommand{\noprint}[1]{}

{\theoremstyle{definition}

\newtheorem*{remark*}{Remark}
}

\newcommand{\checked}[1][\null]{\ensuremath{\boldsymbol{\surd}}}

\begin{document}

\par\noindent {\LARGE\bf
Stochastic domain decomposition\\ for time dependent  adaptive mesh generation
\par}

{\vspace{4mm}\par\noindent {\bf Alexander Bihlo$^\dag\, ^\ddag$, Ronald D.\ Haynes$^\dag$ and Emily J.\ Walsh$^{\S}$
} \par\vspace{2mm}\par}

{\vspace{2mm}\par\noindent {\it
$^{\dag}$~Department of Mathematics and Statistics, Memorial University of Newfoundland,\\
$\phantom{^\dag}$~St.\ John's (NL), A1C 5S7, Canada
}}

{\vspace{2mm}\par\noindent {\it
$^{\ddag}$~Department of Mathematics, University of British Columbia\\ Vancouver (BC), V6T 1Z2, Canada
}}

{\vspace{2mm}\par\noindent {\it
$^{\S}$~Department of Engineering Design and Mathematics, University of Western England, Bristol BS16 1QY, United Kingdom
}}

{\vspace{2mm}\par\noindent {\it
$\phantom{^\dag}$~\textup{E-mail}:
$^{\dag}$abihlo@mun.ca, rhaynes@mun.ca,
$^{\S}$emily3.walsh@uwe.ac.uk
}\par}

\vspace{4mm}\par\noindent\hspace*{8mm}\parbox{140mm}{\small
The efficient generation of meshes is an important component
in the numerical solution of problems in physics and engineering.
Of interest are
situations where global mesh quality and a tight coupling to the solution
of the physical partial differential equation (PDE)
is important.  We consider parabolic PDE  mesh generation
and present a method for the construction of adaptive meshes in two spatial dimensions
using stochastic domain decomposition that is suitable for an implementation
in a multi-- or many--core environment.
Methods for mesh generation on periodic domains are also provided.  The mesh generator is coupled to a time dependent physical PDE and the system is evolved
using an alternating solution procedure.
The method uses the stochastic representation of the exact solution of a
parabolic linear mesh generator  to find the location of an adaptive mesh
along the (artificial) subdomain interfaces.
The deterministic evaluation of the mesh over each subdomain
can then be obtained completely
independently using the probabilistically computed solutions
as boundary conditions.
The parallel performance of this general stochastic domain decomposition approach
has previously been shown.
We demonstrate the approach numerically for the mesh generation context and compare the mesh obtained
with the corresponding single domain mesh using a representative mesh quality
measure.
}\par\vspace{2mm}

\section{Introduction}

The numerical solution of many partial differential equations (PDEs) benefits from the
construction of an adaptive grid automatically tuned by the solution itself.
The quasi--Lagrangian (QL), $r$--refinement, approach used here keeps the number of mesh
points and the mesh topology fixed, moving the mesh continuously in time using a
moving mesh PDE (MMPDE).   The solution of the mesh PDE gives a continuous mesh
transformation between an underlying computational co-ordinate and
the required physical co-ordinate.  Both the mesh and the solution are obtained at each
time step.
The QL approach can be implemented in either an alternating or simultaneous manner.
The simultaneous QL approach treats the MMPDE and physical PDE as one large coupled system.  At each time
the new mesh and new solution on that mesh are found concurrently.
Hence the mesh reacts instantly to changes in the physical solution.
This highly nonlinear coupling may destroy exploitable structure which exists
in the discretization of the physical PDE alone.  The alternating approach uses the current mesh and physical solution to update the mesh alone, this new mesh then facilitates the computation of the updated physical solution.   This
introduces a time lag in the mesh as the new mesh is based only on the current physical solution.
Computationally, however, this decoupling reduces the size of the discrete problem.
Furthermore, the solver becomes more modular; the mesh and physical solvers can be called in alternating fashion;
each solver designed to take advantage of the structure inherent in each subsystem.
The simultaneous approach is generally thought to be more difficult and expensive to solve and
hence the alternating method (or a variant thereof) is typically used in two or more spatial dimensions.
As we will see below, the alternating approach fits well with the stochastic domain decomposition approach we describe
to parallelize our computations.

 The general QL approach
has shown great promise in recent years, solving problems in
 meteorology \cite{budd2012},
relativistic magnetohydrodynamics \cite{he2012},
 combustion and convection in a porous medium \cite{cao99c},
groundwater flow
 and transport of nonaqueous phase liquids
\cite{Huang2005},
Stefan problems \cite{beckett2001},
semiconductor devices \cite{yuan2012},
and viscoelastic flows \cite{zhang2011},
phase change problems \cite{baines2009},
multiphase flows \cite{quan2011}, and low speed viscous flow \cite{jin2010}, to name just
a few.
A thorough overview of PDE based moving mesh methods
 may be found in \cite{huan10a}.

 Recently, motivated by the alternating solution method, one of the authors has studied the parallel solution of the
 nonlinear MMPDE alone using a Schwarz based domain decomposition approach.
 In \cite{Haynes:2012},
Haynes and Gander propose and analyze classical, optimal and optimized
Schwarz methods in one spatial dimension at the continuous level.
A numerical study of
classical and optimized Schwarz domain decomposition for 2D nonlinear mesh generation
has been presented in \cite{hayn12a}.
In \cite{Chen:2012uo}, a  {\em monolithic} domain decomposition method, simultaneously solving a linear
mesh generator coupled to the physical PDE, was presented for a shape
optimization problem.  The authors used an overlapping domain decomposition approach
to solve the coupled problem.

In this paper,  we present an efficient, parallel strategy for the solution of the moving mesh PDE
 based on a stochastic domain decomposition method  proposed by Acebr{\'o}n et al.\ \cite{aceb05a}.
 The motivation is two--fold.  First, we wish to reduce (by parallelization) the potential
 burden of having
 to solve an additional (mesh) PDE.   Second,
it is often mesh quality, not an extremely accurate solution of the mesh PDE, which is
important.  As we will see, the stochastic domain decomposition approach is a means to these ends.

A stochastic domain decomposition (SDD) method to
find adaptive meshes for steady state elliptic problems
was presented in \cite{haynesbihlo2}.
The SDD approach, as originally formulated in~\cite{aceb05a},  uses a
probabilistic form of the point--wise solution for linear elliptic boundary value problems.
The point--wise solution is evaluated only at the
introduced subdomain interfaces.
The approximation of the solution at each interface point is obtained independently using Monte--Carlo simulations and
these evaluations are then
used as Dirichlet
boundary conditions for the (deterministic) subdomain solves which can be computed in parallel.
The mesh PDEs are generally not solved to high accuracy.
Mesh quality, allowing an accurate representation of the physical PDE, is what is generally required.
The parallel algorithm and lower accuracy requirement makes the proposed SDD method computationally attractive.
The lower accuracy requirement allows one to terminate the Monte--Carlo simulations well before convergence.

The existence of a stochastic representation of the
 exact solution of  linear parabolic problems allowed us to extend the
 SDD approach to (linear) parabolic mesh generators in \cite{haynesbihlo2}.
There we considered the time relaxed form of the Winslow--Crowley variable
diffusion mesh generation method, first described  in~\cite{wins66a}.
Only the solution of the mesh generator for a specified analytic mesh density function was considered.

Here we consider the generation of time dependent meshes  where the mesh PDE is coupled to a
physical PDE of interest.
The coupling to the physical solution $u$, is provided by, for example,  an arc-length type
function $\rho=\sqrt{1+\alpha(u_x^2+u_y^2)}$.
Furthermore, due to a stochastic representation for the solution of PDEs subject
to periodic boundary conditions we give, for the first time, a method to
generate meshes for periodic problems using a stochastic domain
decomposition approach.

In Section \ref{sec:meshstuff} we describe the time dependent mesh generator
used in this paper and how the coupling between the physical PDE and
mesh PDE is handled for the global, single domain solution.    Furthermore, we describe a mesh generator for a spatially periodic problem.   Section \ref{sec:stochasticdd} provides some background on the stochastic solution of
linear time dependent, periodic and non--periodic PDEs and describes how this is used to generate a non--iterative domain decomposition algorithm.  We
precisely illustrate how to obtain approximations to  the point--wise solution of the mesh PDEs using the stochastic approach and finally how to generate the meshes using the stochastic domain decomposition framework.
In Section \ref{sec:numerics} we illustrate the algorithm for various examples including Burgers' equation in both the non--periodic and periodic situations and the shallow water equations on a periodic domain.  We conclude with some
observations and items for further study in Section \ref{sec:conclude}.

\section{Mesh generation approach}\label{sec:meshstuff}

As discussed extensively in the references above, in 1D the QL $r$--refinement approach
generates a physical mesh $x\in [0,1]$,  via a continuous mesh transformation $x(\xi)$, where
$\xi$ is an underlying computational variable.
A guiding principle is the equidistribution principle of DeBoor (in 1D) \cite{deboor,burchard,white},
which finds a mesh transformation $x(\xi)$ by enforcing
\[
    \int_0^\xi \rho(t,u,\tilde{x})\,d\tilde{x} = \xi\int_0^1 \rho(t,u,\tilde{x})\,d\tilde{x},
\]
where the mesh density function $\rho$ gives a measure of the level of difficulty or error in the physical
solution $u$.  The physical solution $u$ may be
given analytically or as the solution of a physical
PDE. In differential form, the mesh transformation may be found as the solution of the quasilinear BVP
\begin{equation}\label{eq:nonlinearbvp}
    \frac{d}{d\xi}\left(\rho(x)\frac{dx}{d\xi}\right)=0, \quad x(0)=a,\, x(1)=b.
 \end{equation}
 Here we have suppressed the $t$ and $u$ dependency in $\rho$.
This gives the physical co--ordinates $x_i=x(\xi_i)$
as a function of a (typically uniform) computational grid $\xi_i$.
This BVP can be written as a linear BVP in terms of the inverse mesh transformation $\xi(x)$ as
\begin{equation}\label{eq:linearbvp}
\frac{d}{dx}\left(\frac{1}{\rho(x)}\frac{d\xi}{dx}\right) = 0, \quad \xi(a) = 0,\, \xi(b) = 1.
\end{equation}

The linearity of the BVP for $\xi(x)$ makes it easier to solve (in some sense)
and in higher dimensions has the additional benefit that it is easy to say concrete things about
the well--posedness of the mesh transformation. As will be discussed below, the linearity of a mesh generator is also a prerequisite for the stochastic domain decomposition method. The obvious disadvantage is that the solution of the
BVP for $\xi(x)$ does not give the mesh locations in the physical variables $x$.  One alternative
is to transform the physical PDE from the physical variables to the new computational co-oordinates.
Alternatively, inverse linear interpolation could be used to find the $x$ locations by projecting
$\xi(x)$ onto a uniform $\xi$ grid.

Indeed, in \cite{haynesbihlo2} we constructed our stochastic DD method for steady state mesh generation using the natural 2D extension of the linear BVP (\ref{eq:linearbvp}).  In the time dependent case considered in this paper we will construct a stochastic DD method directly in the physical co-ordinates.  The time stepping will provide a natural linearization as we will show below.

A natural way to derive the BVPs above (and which extends to higher dimensions) is as the Euler--Lagrange equations whose solutions
minimize certain functionals of the required mesh transformations.
For example minimizing the functional
\[
I[x] = \frac{1}{2}\int_0^1 \left(\rho(x)\frac{dx}{d\xi}\right)^2\,d\xi
\]
leads to the quasi--linear BVP (\ref{eq:nonlinearbvp}) above.   A minimizer $\xi(x)$ of the functional
\[
I[\xi] = \frac{1}{2}\int_a^b \frac{1}{\rho(x)}\left(\frac{d\xi}{dx}\right)^2\,dx
\]
satisfies the linear BVP (\ref{eq:linearbvp}) above.



For time dependent PDEs it is useful to use
a formulation which involves the mesh speed, such a mesh equation is called a moving mesh PDE (MMPDE).
If $\rho = \rho(t,x)$ the functional $I[\xi]$ becomes
\[
I[\xi] = \frac{1}{2}\int_a^b \frac{1}{\rho(t,x)}\left(\frac{d\xi}{dx}\right)^2\,dx.
\]
As described in \cite{huan10a}, the direction $\xi$ which reduces $I[\xi]$ is given by the gradient flow equation
\[
\frac{d\xi}{dt} = \frac{P}{\tau}\frac{d}{dx}\left(\frac{1}{\rho(t,x)}\frac{d\xi}{dx}\right),
\]
where $\tau > 0$ is a user specified constant which determines the response of the mesh to changes in
$\rho$ (or $u$).  Here $P$ is a positive definite differential operator that we choose with some
flexibility.


We can switch the roles
of $x$ and $\xi$ in the gradient flow equation to get a mathematically equivalent moving mesh PDE in the
physical variables $x$:
\[
\frac{\partial x}{\partial t} = \frac{1}{\tau}\frac{\partial x}{\partial \xi}P\left(\rho\frac{\partial x}{\partial \xi}\right)^{-2}\left(\frac{\partial x}{\partial \xi}\right)^{-1}\frac{\partial}{\partial \xi}\left(\rho\frac{\partial x}{\partial \xi}\right).
\]
Hence this resulting equation would retain the nice well--posedness properties.

If we choose $P = (\rho x_\xi)^2$, then we get
\begin{equation}\label{eq:mmpde5}
\frac{\partial x}{\partial t} = \frac{1}{\tau}\frac{\partial}{\partial\xi}\left(\rho\frac{\partial x}{\partial \xi}\right),
\end{equation}
which is referred to as MMPDE5.
Our focus here is the generation of time dependent meshes by using this
nonlinear parabolic mesh generator.

The coupling to the physical PDE and physical solution $u$
is provided by the mesh density function $\rho(x,u)$.
In a typical deterministic implementation, the mesh  and physical solution are
updated by discretizing
MMPDE5 and the physical PDE in time and either solving the resulting large
nonlinear system of algebraic equations for both the mesh and physical solution
simultaneously or proceeding in an alternating fashion.   The alternating or MP
approach \cite{huan10a} used in this paper freezes
$\rho$ in the time discretized mesh equation
at the current $u^n$ to compute the next mesh $x^{n+1}$ and then
integrates the physical PDE, using mesh $x^{n+1}$, to obtain $u^{n+1}$.
The stochastic solution representation given below requires the PDE to be
linear.   This MP procedure effectively linearizes the MMPDE.
In two spatial dimensions we use  the mesh generator
\begin{equation}\label{eq:ParabolicWinslowGeneratorPeriodicDomain}
\begin{split}
 &x_t=\frac{\nabla_{\boldsymbol\xi}\rho}{\rho}\cdot\nabla_{\boldsymbol\xi}x+\nabla_{\boldsymbol \xi}^2x,\quad y_t=\frac{\nabla_{\boldsymbol\xi}\rho}{\rho}\cdot\nabla_{\boldsymbol\xi}y+\nabla_{\boldsymbol \xi}^2y,\\
 &x(0,\boldsymbol \xi)=x_0=L_x\xi,\quad y(0,\boldsymbol \xi)=y_0=L_y\eta,
\end{split}
\end{equation}
where $\boldsymbol{\xi}=(\xi, \eta)$ and $\nabla_{\boldsymbol \xi}=(\partial /\partial \xi, \partial /\partial \eta)$, which we solve over the
square computational domain $\Omega_{\rm c}=[0,1]\times[0,1]$.
For the sake of simplicity we assume a rectangular domain $\Omega_{\rm p}=[0,L_x]\times [0,L_y]$, where $L_x>0$, $L_y>0$.
At the actual boundaries of the computational domain $\Omega_{\rm c}$, we employ the fixed boundary conditions $x(0,\eta)=0$, $x(1,\eta)=L_x$, $y(\xi,0)=0$ and $y(\xi,1)=L_y$. The remaining boundary conditions, $x(\xi,0)$, $x(\xi,1)$, $y(0,\eta)$
and $y(1,\eta)$ are found by solving the respective one-dimensional forms of the mesh generator~(\ref{eq:ParabolicWinslowGeneratorPeriodicDomain})  along these boundaries.

Note that each of $x$ and $y$ in  mesh generator~(\ref{eq:ParabolicWinslowGeneratorPeriodicDomain})  may be obtained by \emph{first} dividing the two dimensional version of ~(\ref{eq:nonlinearbvp}) by~$\rho$ and \emph{then} relaxing the equation in time. We have found in practice that this form of a relaxed mesh generator gives better meshes than by relaxing the original version~(\ref{eq:nonlinearbvp}). In the discrete formulation, these two possible time relaxed mesh generators are related by a scaled time step.


We are also interested in mesh generation on periodic domains.
It has been pointed out in~\cite{tan08a}, that in order to use the mesh generator~\eqref{eq:ParabolicWinslowGeneratorPeriodicDomain} on a rectangular periodic domain of physical dimensions $\Omega_{\rm p}=[0,L_x[\times[0,L_y[$,
 \[
  x(t,1,\eta)=x(t,0,\eta)+L_x,\quad y(t,\xi,1)=y(t,\xi,0)+L_y,
 \]
 one should express the mesh and physical PDEs in terms of the displacements $\mathbf{X}=\mathbf{x}-\boldsymbol{\bar{\xi}}$, where $\boldsymbol{\bar{\xi}}=(L_x\xi, L_y\eta)$, $\mathbf{x}=(x,y)$ and $\mathbf{X}=(X,Y)$, which are directly periodic on $\Omega_{\rm c}$. In this new set of variables, the mesh generator~\eqref{eq:ParabolicWinslowGeneratorPeriodicDomain} becomes
\begin{align*}
  &X_t=\frac{\nabla_{\boldsymbol\xi}\rho}{\rho}\cdot\nabla_{\boldsymbol\xi}X+\nabla_{\boldsymbol \xi}^2X+\frac{\rho_\xi}{\rho},\quad Y_t=\frac{\nabla_{\boldsymbol\xi}\rho}{\rho}\cdot\nabla_{\boldsymbol\xi}Y+\nabla_{\boldsymbol \xi}^2Y+\frac{\rho_\eta}{\rho}.
\end{align*}

An alternative to the procedure proposed in~\cite{tan08a} works directly in the original variables $x$ and $y$ and properly extends the computational domain using the periodicity of the grid. This allows one to evaluate the derivatives of $\mathbf{x}$ on the boundaries without having to take into account the actual size of the physical domain $\Omega_{\rm p}$. As this approach allows one to work with the physical coordinates $\mathbf{x}$ directly, we
use this second method in this paper.

\color{black}

\section{Stochastic domain decomposition and mesh generation}
\label{sec:stochasticdd}

In this section we describe stochastic domain decomposition methods
for mesh generation on periodic and non--periodic domains and describe an implementation which couples
the mesh generator to the solution of a physical PDE.

\subsection{Stochastic domain decomposition for linear PDEs}

The main idea of stochastic domain decomposition as proposed in~\cite{aceb05a} (see also~\cite{aceb07a}) rests on the stochastic representation of the exact solution to  \emph{linear elliptic} (resp.\ \emph{parabolic}) \emph{boundary value problems}. This now classical connection between stochastic analysis and boundary value problems was first uncovered by Kakutani~\cite{kaku44a,kaku44b} for the Dirichlet problem using Brownian motion. The monographs \cite{kara91a,okse10a} provide
a more recent exposition.

Numerically evaluating this stochastic representation of the exact solution of a linear boundary value problem using Monte--Carlo methods enables one to compute the point-wise numerical solution to the underlying PDE. From the practical point of view, this is fundamentally different to solving a PDE using, say, finite differences, which requires the computation of the numerical solution over the entire domain even if it is only needed in a single point.

Notoriously, Monte--Carlo methods converge slowly, with convergence rates proportional to $N^{-1/2}$ where
$N$ is the number of Monte--Carlo simulations, if pseudo-random numbers are used~\cite{pres07a}, and hence they play a role mostly for higher dimensional problems, where they can be shown to outperform deterministic methods.

An alternative is to use them in the context of domain decomposition. Namely, splitting the entire domain into non-overlapping subdomains, the stochastic solution can be used to compute the point-wise interface solutions between the subdomains. Once these interface solutions are determined with sufficient accuracy, they act as Dirichlet boundary values for the single subdomains. The PDE solution over each subdomain is computed deterministically using an appropriate discretization of the underlying PDE. The main advantage of this method is that iteration, as is required in classical domain decomposition methods (such as Schwarz methods), can be completely avoided. Also, the solutions over each subdomain can be obtained in parallel and thus the method is suitable for massively parallel computing architectures.

In~\cite{haynesbihlo2,bihl15a} we have shown that the stochastic domain decomposition technique is an effective way for the parallel generation of adaptive meshes. A main motivator for the approach is that it is, in general, not necessary to compute the meshes with high accuracy. What is important is to obtain meshes with high mesh quality. It was shown in~\cite{haynesbihlo2,bihl15a} that even meshes that are not accurate solutions to the mesh PDEs
can have good mesh quality. This characteristic of mesh generation enables an increase in the efficiency of the stochastic domain decomposition method.

\subsection{Stochastic analysis for Dirichlet boundary value problems}\label{sec:SDDDirichletBVP}

In this section we give the specifics of the stochastic analysis required to
generate the solution of linear parabolic boundary value problems.  Specifically, we consider system (\ref{eq:ParabolicWinslowGeneratorPeriodicDomain})
where  the unknowns $x=x(t,\xi,\eta)$ and $y=y(t,\xi,\eta)$ are required on the time-space domain given by $[0, T]\times\Omega_{\rm c}$, with $T$ being some finite final time.
System~\eqref{eq:ParabolicWinslowGeneratorPeriodicDomain} is supplemented with the boundary conditions $x(t,\xi,\eta)|_{\partial \Omega_{\rm c}}=f(t,\xi,\eta)$ and $y(t,\xi,\eta)|_{\partial \Omega_{\rm c}}=g(t,\xi,\eta)$, for given continuous functions $f$ and $g$.  The initial values are $x(0,\xi,\eta)=x_0(\xi,\eta)$ and $y(0,\xi,\eta)=y_0(\xi,\eta)$.

It is important to stress here that \eqref{eq:ParabolicWinslowGeneratorPeriodicDomain} is nonlinear if the mesh density function is a function on the physical domain, that is
$\rho=\rho(t,x,y)$. However, as was indicated in Section~\ref{sec:meshstuff}, in the practical implementation it is possible to freeze~$\rho$ at time layer $t^n$ when computing the mesh at time $t^{n+1}$. This effectively boils down to a linearization of the mesh generator~\eqref{eq:ParabolicWinslowGeneratorPeriodicDomain}. It is thus appropriate to assume that $\rho$ in~\eqref{eq:ParabolicWinslowGeneratorPeriodicDomain} is not a function of~$x$ and~$y$ but of some auxiliary variables~$\tilde x$ and~$\tilde y$ (and time), i.e.\ we assume that $\rho=\rho(t,\tilde x,\tilde y)$. Then this mesh generator becomes linear and allows for a stochastic representation of its exact solution \cite{mils04a}, given by
\begin{subequations}\label{eq:WinslowGeneratorRelaxedStochasticSolution}
\begin{align}\label{eq:WinslowGeneratorRelaxedStochasticSolutionA}
\begin{split}
 &x(t,\xi,\eta)=\mathrm{E}\left[x_0(\mathbf{\Phi}(t))\mathbf{1}_{[\tau_{\partial\Omega_{\rm c}}>t]} \right] + \mathrm{E}\left[f(t-\tau_{\partial\Omega_{\rm c}}, \mathbf{\Phi}(\tau_{\partial\Omega_{\rm c}}))\mathbf{1}_{[\tau_{\partial\Omega_{\rm c}}<t]} \right],\\
 &y(t,\xi,\eta)=\mathrm{E}\left[y_0(\mathbf{\Phi}(t))\mathbf{1}_{[\tau_{\partial\Omega_{\rm c}}>t]} \right] + \mathrm{E}\left[g(t-\tau_{\partial\Omega_{\rm c}}, \mathbf{\Phi}(\tau_{\partial\Omega_{\rm c}}))\mathbf{1}_{[\tau_{\partial\Omega_{\rm c}}<t]} \right],
\end{split}
\end{align}
where the stochastic process $\mathbf{\Phi}(t)$ satisfies the stochastic differential equation (SDE)
\begin{equation}\label{eq:WinslowGeneratorRelaxedStochasticSolutionB}
    \mathrm{d}\mathbf{\Phi}(t)=\frac{1}{\rho}\nabla_{\boldsymbol\xi} \rho\,\mathrm{d}t+\sqrt{2}\,\mathrm{d}\mathbf{W}(t).
\end{equation}
\end{subequations}
In~\eqref{eq:WinslowGeneratorRelaxedStochasticSolution}, $\mathrm{E}[\cdot]$ is the expected value, $\tau_{\partial\Omega_{\rm c}}$ is the first exit time of the stochastic process $\mathbf{\Phi}(t)$ starting at $(\xi,\eta)$, $\mathbf{W}$ is two-dimensional Brownian motion and $\mathbf{1}$ is the indicator function. See~\cite{kara91a,mils04a,okse10a} for a more extensive discussion of this subject.
We note that the stochastic solution~\eqref{eq:WinslowGeneratorRelaxedStochasticSolutionA} has contributions from both the specific initial condition and boundary values of the mesh generator.

\subsection{Stochastic analysis for periodic problems}\label{sec:GlobalPeriodicProblem}

The mesh generator~\eqref{eq:ParabolicWinslowGeneratorPeriodicDomain} is in the form of a system of linear, second order, parabolic PDEs. In Section~\ref{sec:SDDDirichletBVP} we wrote the stochastic point--wise solution assuming
a Dirichlet boundary value problem.
On periodic domains, the celebrated Kac--Feynman formula can be used to obtain the stochastic representation of the mesh generator~\eqref{eq:ParabolicWinslowGeneratorPeriodicDomain}, see e.g.~\cite{mils04a}. In this case, only initial values $x(0,\xi,\eta)=x_0(\xi,\eta)$ and $y(0,\xi,\eta)=y_0(\xi,\eta)$ are given. The solution to~\eqref{eq:ParabolicWinslowGeneratorPeriodicDomain} can then be written as
\begin{equation}\label{eq:StochasticSolutionPeriodicDomainSolution}
 x(t,\xi,\eta)=\mathrm{E}\left[x_0(\boldsymbol{\Phi})\right],\quad
 y(t,\xi,\eta)=\mathrm{E}\left[y_0(\boldsymbol{\Phi})\right],
\end{equation}
where $\boldsymbol{\Phi}$ satisfies the same stochastic differential equation as given in~\eqref{eq:WinslowGeneratorRelaxedStochasticSolutionB}.

\subsection{Stochastic domain decomposition on periodic and non--periodic domains}

While the stochastic solutions~\eqref{eq:WinslowGeneratorRelaxedStochasticSolution} and \eqref{eq:StochasticSolutionPeriodicDomainSolution} can in principle be used to obtain the solution to the parabolic mesh generator at any time $t>0$, our mesh generation problem is slightly more complicated.
The mesh density function $\rho$ is linked to the solution of the physical PDE system, which changes over time. For this reason, in practice we use~\eqref{eq:WinslowGeneratorRelaxedStochasticSolution} and \eqref{eq:StochasticSolutionPeriodicDomainSolution} only to advance the mesh over a single time step, $\Delta t$,  from $t^n$ to $t^{n+1}$.

We discretize the SDE~\eqref{eq:WinslowGeneratorRelaxedStochasticSolutionB} using the Euler-Maruyama method,
\begin{equation}\label{eq:WinslowGeneratorSDEDiscretization}
 \boldsymbol{\Phi}^{k+1}=\boldsymbol{\Phi}^{k}+\frac{\nabla_{\boldsymbol\xi}\rho}{\rho}\bigg|_{(t^n,\boldsymbol{\Phi}^k)}\Delta t_{\rm s} + \sqrt{2\Delta t_{\rm s}}\mathbf{N}(0,1),
\end{equation}
with constant time step $\Delta t_{\rm s}=\Delta t/K$, $k=0,\dots,K-1$, where $\mathbf{N}$ is a two-dimensional vector of Gaussian distributed random numbers with zero mean and variance one. In other words, one time step of size $\Delta t$ is split into $K$ sub-time steps. This splitting is necessary so that the use of an excessively small time-step $\Delta t$ for the solution of the physical differential equation can be avoided. The mesh density function $\rho$ remains fixed at time step $t^n$.  The derivatives $\nabla_{\boldsymbol{\xi}}\rho$ are approximated with finite differences.

In practice, the starting points of the stochastic process at time $t^n$, $\boldsymbol{\Phi}^0$, coincide with the grid points where the stochastic solution is required. Once the new values $\boldsymbol{\Phi}^{k+1}$ at time $t^{k+1}$ are computed, both $\rho$ and $\nabla_{\boldsymbol{\xi}}\rho$ have to be approximated at $\boldsymbol{\Phi}^{k+1}$. For this, bi-linear interpolation is used. The procedure is repeated until $\boldsymbol{\Phi}^{K}$ is computed, which coincides with the value of the stochastic process at time $t^{n+1}$. We then evaluate the initial values of $x_0$ and $y_0$ given at time $t^n$ at the new location $\boldsymbol{\Phi}^K$ using bi-cubic interpolation. This gives the values $\tilde x_0(\boldsymbol{\Phi}^K)$ and $\tilde y_0(\boldsymbol{\Phi}^K)$, which approximate the actual values $x_0(\boldsymbol{\Phi}^{n+1})$ and $y_0(\boldsymbol{\Phi}^{n+1})$ that are needed in both the solutions~\eqref{eq:WinslowGeneratorRelaxedStochasticSolutionA} (for the first term) and~\eqref{eq:StochasticSolutionPeriodicDomainSolution}.

Solving  Dirichlet boundary value problems stochastically, a boundary test has to be applied to determine whether the process $\boldsymbol{\Phi}^{k+1}$ has left the domain within the sub-time step $[t^k,t^{k+1}]$. A linear Brownian bridge is used as interpolating process as discussed in~\cite{jans03a}. If the process left the domain before time $t^{n+1}$, the integration can be stopped and the second term in \eqref{eq:WinslowGeneratorRelaxedStochasticSolutionA}  can be evaluated. No such boundary test is needed for periodic domains.

In order to estimate the expected value using the arithmetic mean, the procedure is repeated $N$ times.

\subsection{Local subdomain problem}

In the domain decomposition solution to the problem~\eqref{eq:ParabolicWinslowGeneratorPeriodicDomain}, we only use the stochastic solutions
\eqref{eq:WinslowGeneratorRelaxedStochasticSolution} and \eqref{eq:StochasticSolutionPeriodicDomainSolution} to generate the subdomain interface solutions for the Dirichlet and
periodic problems. Once these solutions are computed, the solution to the mesh generation problem on the individual subdomains becomes a Dirichlet boundary value problem.
In order to solve this problem, the original parabolic mesh generator~\eqref{eq:ParabolicWinslowGeneratorPeriodicDomain} has to be solved with the Dirichlet boundary conditions
\[
 \mathbf{x}|_{\partial \Omega_{\rm c}^i}=f^i,
\]
where $\partial \Omega_{\rm c}^i$ denotes the boundary of the $i$--th subdomain of the
computational domain and $f^i$ are the values obtained either
from the stochastic representation~\eqref{eq:WinslowGeneratorRelaxedStochasticSolution} or the
boundary conditions on the global mesh where $\partial\Omega_{\rm c}^i \bigcap \partial\Omega_c \neq \emptyset.$

For the local subdomain solver we discretize~\eqref{eq:ParabolicWinslowGeneratorPeriodicDomain} using centered finite differences for the spatial derivatives and a forward Euler method for the time stepping.

\subsection{Domain decomposition solution}

The domain decomposition strategy relies on combining the stochastic evaluation along the interfaces with the deterministic subdomain solves. At each time step, the stochastic solution procedure provides the Dirichlet boundary conditions for the deterministic single-domain solver. This allows the computation of the new mesh at the time step $t^{n+1}$.    The solution values can either be evaluated stochastically
at each $(\xi_i,\eta_j)$ along the artificial interfaces or a sample can be evaluated with the rest obtained by interpolation.   An optimal placement strategy for the interpolation nodes was presented in \cite{haynesbihlo2}.

Once the new mesh is computed, the mesh density function $\rho(t^{n+1},x,y)$ is evaluated on the new mesh and the solution procedure is repeated with the new mesh density function.

\subsection{Mesh quality}

As mentioned previously, there is a  crucial difference between the SDD method as
proposed in~\cite{aceb05a} for the computation of numerical solutions of general
linear elliptic boundary value problems and for the mesh generation case.  The latter
does not necessarily require a very accurate solution of the mesh equation; here, mesh quality is more important.

There are several ways of assessing the quality of an adaptive mesh and different mesh quality measures based on properties such as equidistribution and alignment have been derived. See~\cite{huan10a} for an extensive discussion of mesh quality measures. In~\cite{haynesbihlo2} it was proposed to use the geometric mesh quality measure for the assessment of the quality of meshes generated using the SDD method. The geometric mesh quality measure is defined by
\begin{equation}\label{eq:Qmeasure}
 Q(K)=\frac12\frac{\mathrm{tr}(J^{\rm T}J)}{\sqrt{\det(J^{\rm T}J)}},
\end{equation}
where $J$ is the Jacobian of the transformation $x=x(\xi,\eta)$, $y=y(\xi,\eta)$ and $K$ is a mesh element in $\Omega_{\rm c}$. The quantity $Q(K)$
measures how far the mesh cell is from being equilateral,  that is  we have $Q(K)\ge1$ with $Q(K)=1$ precisely for an equilateral cell.

In~\cite{haynesbihlo2} we argued that using this measure is appropriate as meshes computed using stochastic methods show several kinks when they are far away from convergence. Thus, the values of $Q(K)$ are in general larger for such meshes compared to grids that are computed using deterministic methods. Of course, the absolute value of $Q(K)$ depends on the mesh density function used and the underlying problem for which an adaptive grid is
computed. We are therefore not interested in the absolute value of $Q(K)$ but only in the ratio between $Q^{\rm SD}(K)$, the geometric mesh quality measure of the reference mesh obtained on a single domain,
and $Q^{SDD}(K)$, the geometric mesh quality measure of the grid obtained using the SDD method. If this ratio $R(K)=Q^{\rm SD}/Q^{\rm SDD}$ is close to one, the mesh obtained from the SDD method is a good approximation to the single domain mesh.

The quantity $R$ is a  function of each individual mesh cell. In the next section, for the sake of convenience, we only list the maximum and mean values, $R_{\max}$ and $R_{\mathrm{mean}}$ respectively,  over all mesh cells.

\section{Numerical Results}\label{sec:numerics}
In this section, numerical experiments in 2D
are presented to demonstrate the effectiveness of the domain decomposition algorithm and its suitability for problems with both Dirichlet and periodic boundary conditions.

\subsection{Burgers' equation with Dirichlet boundary conditions}

\begin{figure}[!ht]
\begin{subfigure}{.33\textwidth}
  \centering
  \includegraphics[width=4.5cm,height=4.5cm]{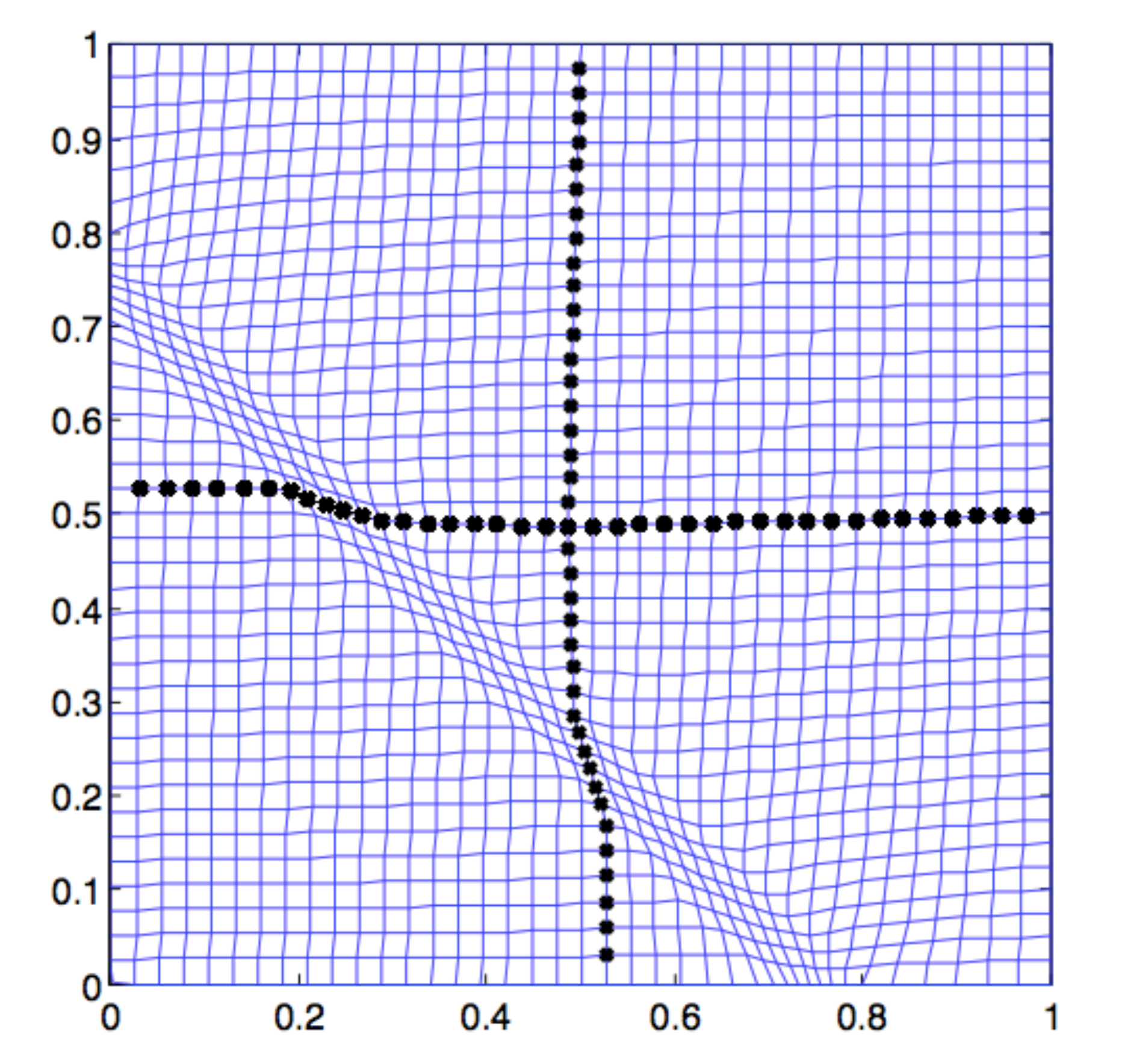}
  \caption{t=0.75}
  \label{fig:Db1}
\end{subfigure}%
\begin{subfigure}{.33\textwidth}
  \centering
  \includegraphics[width=4.5cm,height=4.5cm]{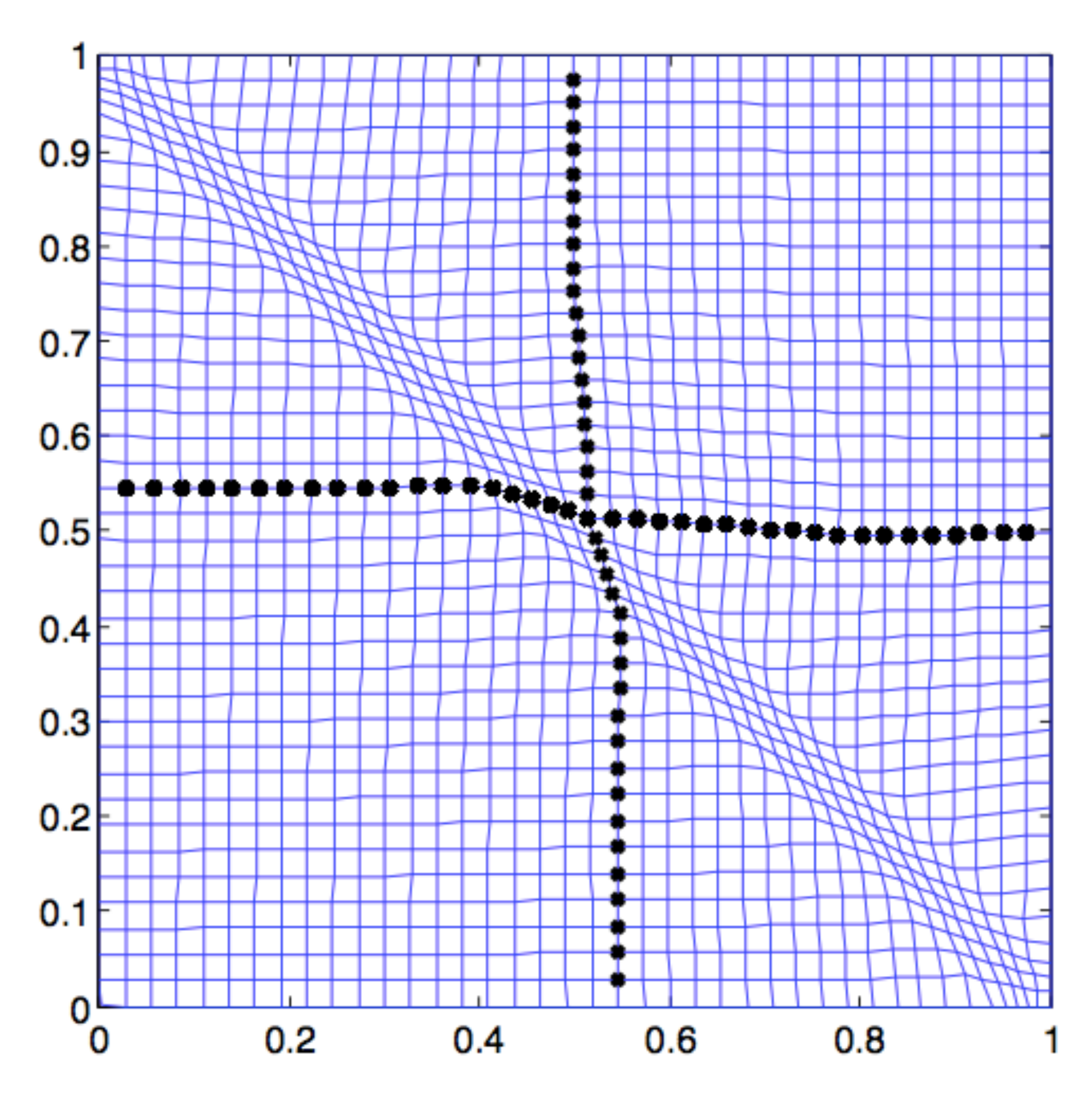}
  \caption{t=1.00}
  \label{fig:Db2}
  \end{subfigure}
 \begin{subfigure}{.33\textwidth}
  \centering
  \includegraphics[width=4.5cm,height=4.5cm]{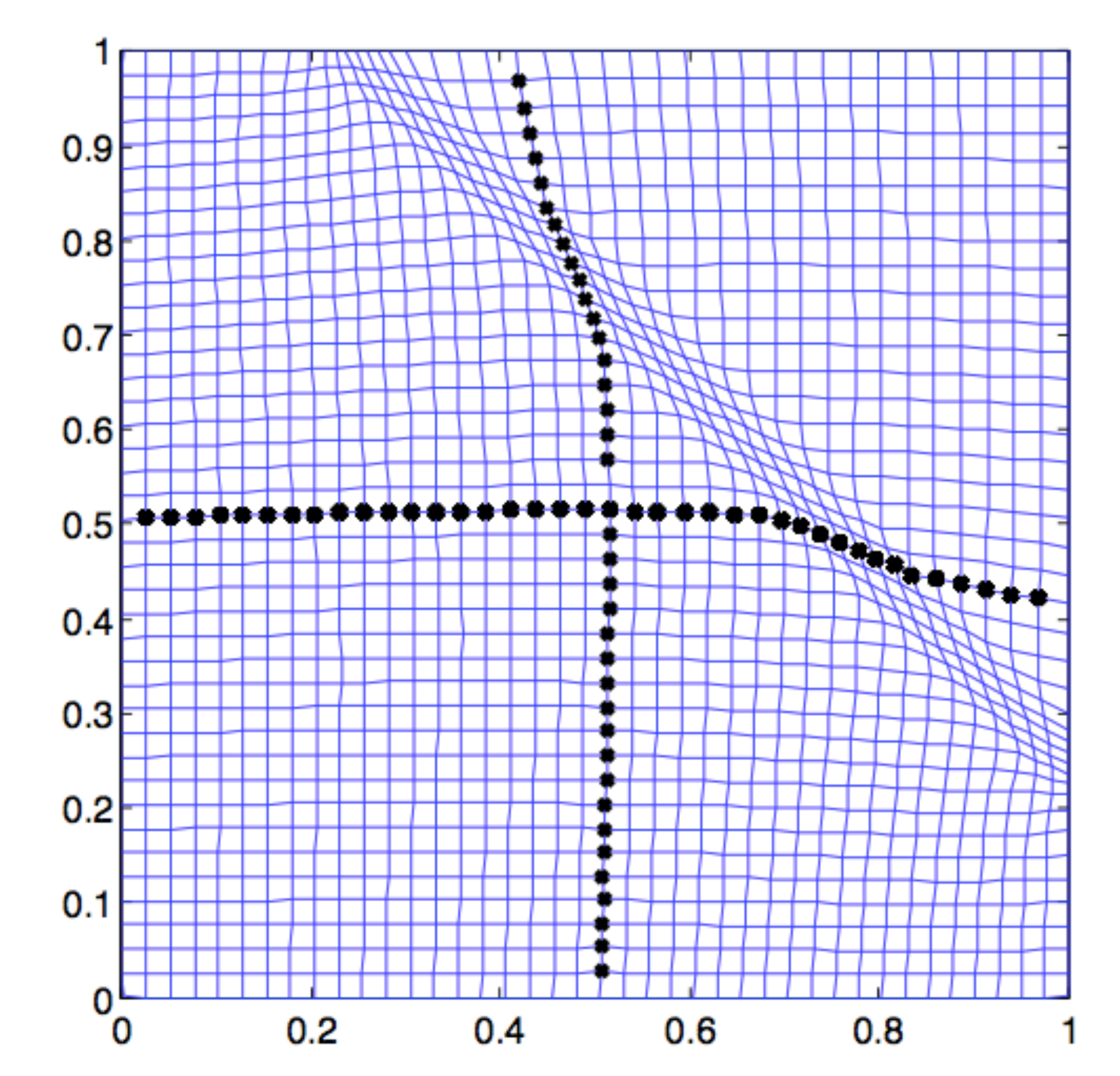}
  \caption{t=1.25}
  \label{fig:Db2}
  \end{subfigure} \\
  \begin{subfigure}{.33\textwidth}
  \centering
  \includegraphics[width=4.5cm,height=4.5cm]{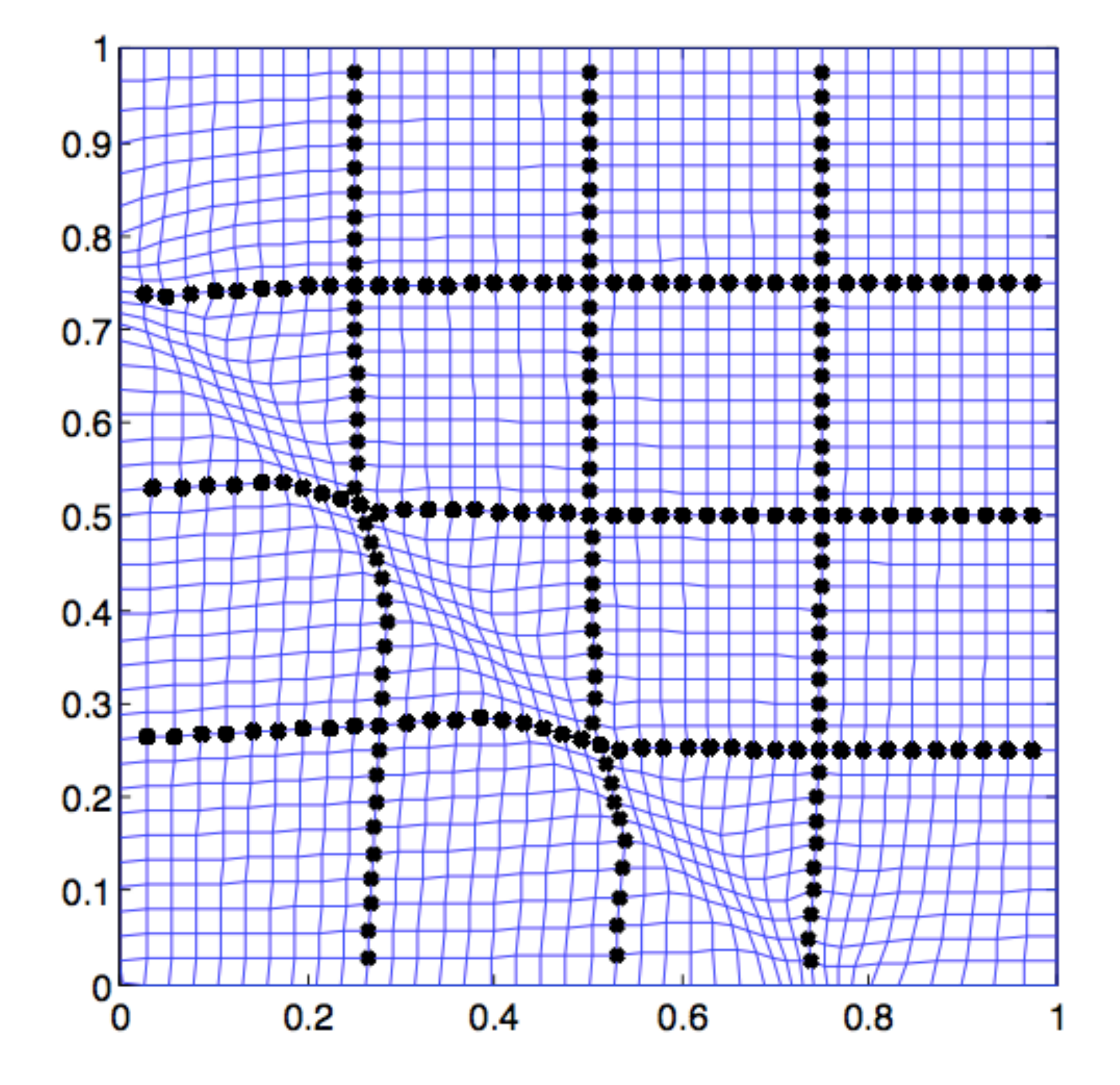}
  \caption{t=0.75}
  \label{fig:Db1}
\end{subfigure}%
\begin{subfigure}{.33\textwidth}
  \centering
  \includegraphics[width=4.5cm,height=4.5cm]{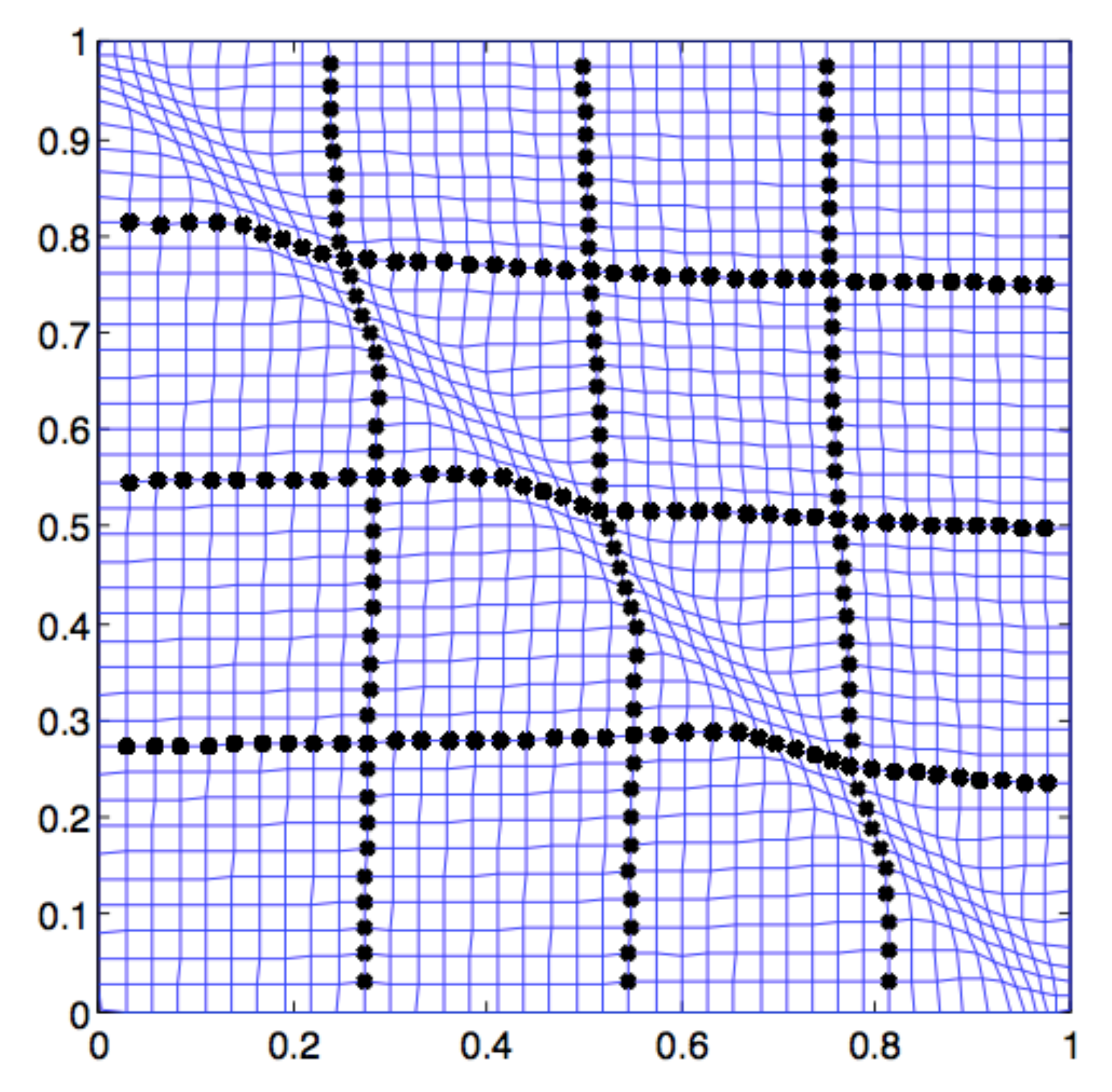}
  \caption{t=1.00}
  \label{fig:Db2}
  \end{subfigure}
 \begin{subfigure}{.33\textwidth}
  \centering
  \includegraphics[width=4.5cm,height=4.5cm]{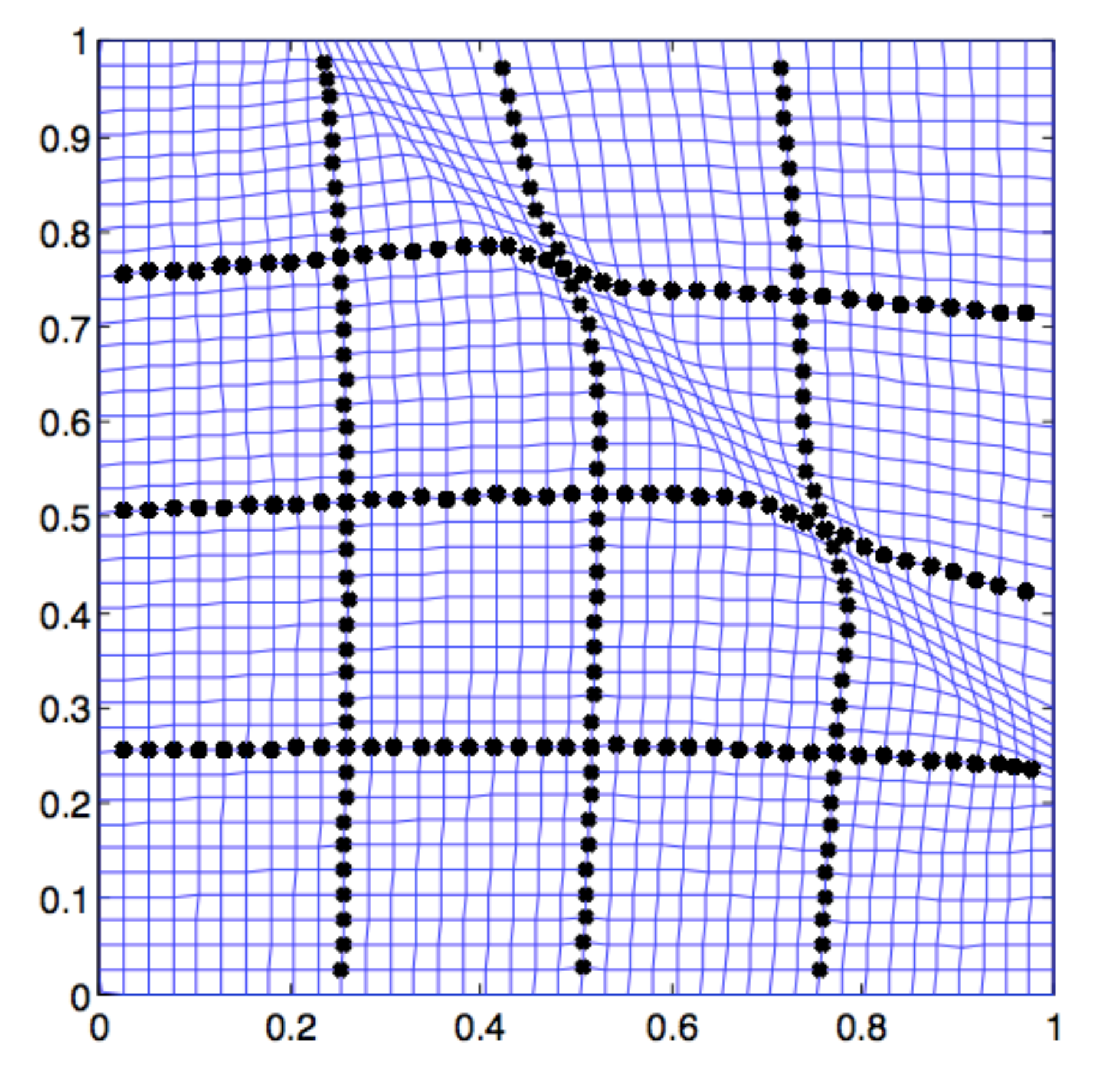}
  \caption{t=1.25}
  \label{fig:Db2}
  \end{subfigure} \\
\caption{Evolution of mesh for Burgers' equation using the stochastic DD method with $2\times2$ (top) and $4\times4$ (bottom) subdomains.}
\label{fig:Dburgers_figs}
\end{figure}

The first test problem considered is the scalar form of the two-dimensional Burgers' equation
\begin{equation}\label{eq:BurgersEquation2D}
 u_t + (u^2/2)_ x+ (u^2/2)_y + \nu (u_{xx}+u_{yy})=0,
\end{equation}
$[x,y]\in\Omega_p=[0,1]\times[0,1]$. The initial and Dirichlet boundary conditions are chosen such that the exact solution is
\[
 u=\left(1+\exp\left(\frac{x+y-t_0}{2\nu}\right)\right)^{-1},
\]
and we consider the case of a moderately small diffusion coefficient $\nu=0.005$. In this problem the smaller $\nu$, the more convection dominates, and large gradients develop and move
to the boundaries for $t > 0$, thus requiring a higher concentration of mesh points to resolve the oblique shock that propagates with time. This is a common test problem in the moving mesh literature \cite{lang,zhang}. The coordinate transformation $(x,y) = (x(\xi,t),y(\xi,t))$ is used  to rewrite the 2D Burgers' equation in QL form in computational coordinates, and this is discretized in the computational domain using centered differences in space and a trapezoidal rule for the time integration. This is coupled
to the mesh generator \eqref{eq:ParabolicWinslowGeneratorPeriodicDomain}. Both (\ref{eq:ParabolicWinslowGeneratorPeriodicDomain}) and \eqref{eq:BurgersEquation2D} are then solved alternately in time (the MP procedure) for each subdomain using an arc-length mesh density function
\[
 \rho=\sqrt{1+\alpha(u_x^2+u_y^2)},
\]
with $\alpha=10/\lVert u_x^2+u_y^2\lVert_{\infty}$. The Dirichlet boundary conditions (at the subdomain interfaces) are found by evaluating
 \eqref{eq:WinslowGeneratorRelaxedStochasticSolutionA} using \eqref{eq:WinslowGeneratorSDEDiscretization}. At the physical boundary a 1D version of
 (\ref{eq:ParabolicWinslowGeneratorPeriodicDomain}) is solved.
We integrate the finite difference discretization of~\eqref{eq:BurgersEquation2D} using $41\times 41$ grid points on the physical domain $\Omega_{\rm c}=[0,1]\times[0,1]$ with a time step $\Delta t=0.001$. The initial time is $t_0=0.25$. We use $N= 10 000$ Monte--Carlo simulations with a time step $\Delta t_s=\Delta t/10$ for the solution of the stochastic differential equation~\eqref{eq:WinslowGeneratorSDEDiscretization}. In Fig.{\ref{fig:Dburgers_figs}} the evolution of the mesh is shown when the physical domain is split into $2\times 2$ and $4\times 4$ subdomains.

We report the mesh qualities and the $l_\infty$-norm comparing the exact solution to the numerical solution at times $t_f=0.75$, $t_f=1$ and $t_f=1.25$ in Table~\ref{tab:BurgersDirichlet}. The geometric mesh quality measure is computed for the single domain solution and compared to the DD solution obtained by splitting the physical domain into $2\times 2$, $3\times 3$ and $4\times 4$ subdomains.

\begin{table}[!ht]
\centering
\caption{Mesh quality and $l_\infty$-errors for Burgers' equation with Dirichlet boundary conditions.}
\begin{tabular}{|r|r||r|r|r||r|r|r||r|r|r|}
  \hline
  $t_f$ & $l^{\rm SD}_\infty$ & $l^{2\times2}_\infty$ & $R^{2\times2}_{\max}$ & $R^{2\times2}_{\mathrm{mean}}$ & $l^{3\times3}_\infty$ & $R^{3\times3}_{\max}$ & $R^{3\times3}_{\mathrm{mean}}$ & $l^{4\times4}_\infty$ & $R^{4\times4}_{\max}$ & $R^{4\times4}_{\mathrm{mean}}$ \\
  \hline\hline
  $0.75$ & 0.027 & 0.023 & 0.99 & 0.99 & 0.031 & 0.95 & 0.98 & 0.026 & 0.95 & 0.98 \\
  $1$    & 0.031 & 0.033 & 0.99 & 0.99 & 0.320 & 0.97 & 0.99 & 0.034 & 0.97 & 0.99 \\
  $1.25$ & 0.030 & 0.031 & 1 & 1 & 0.35 & 0.92 & 1 & 0.24 & 1 & 1 \\
  \hline
\end{tabular}
\label{tab:BurgersDirichlet}
\end{table}

Table~\ref{tab:BurgersDirichlet} shows that all the $l_\infty$-errors
on meshes generated using SDD  are approximately equal to the errors found on the associated
single domain meshes. Also, the ratios of the geometric mesh quality measures of the single domain and DD solutions are close to one for all times, indicating that the DD solutions are a good approximation to the single domain solution.

\subsection{Periodic Mesh Generation}

\begin{figure}[!ht]
\begin{subfigure}{.5\textwidth}
  \centering
  \includegraphics[width=4.5cm,height=4.5cm]{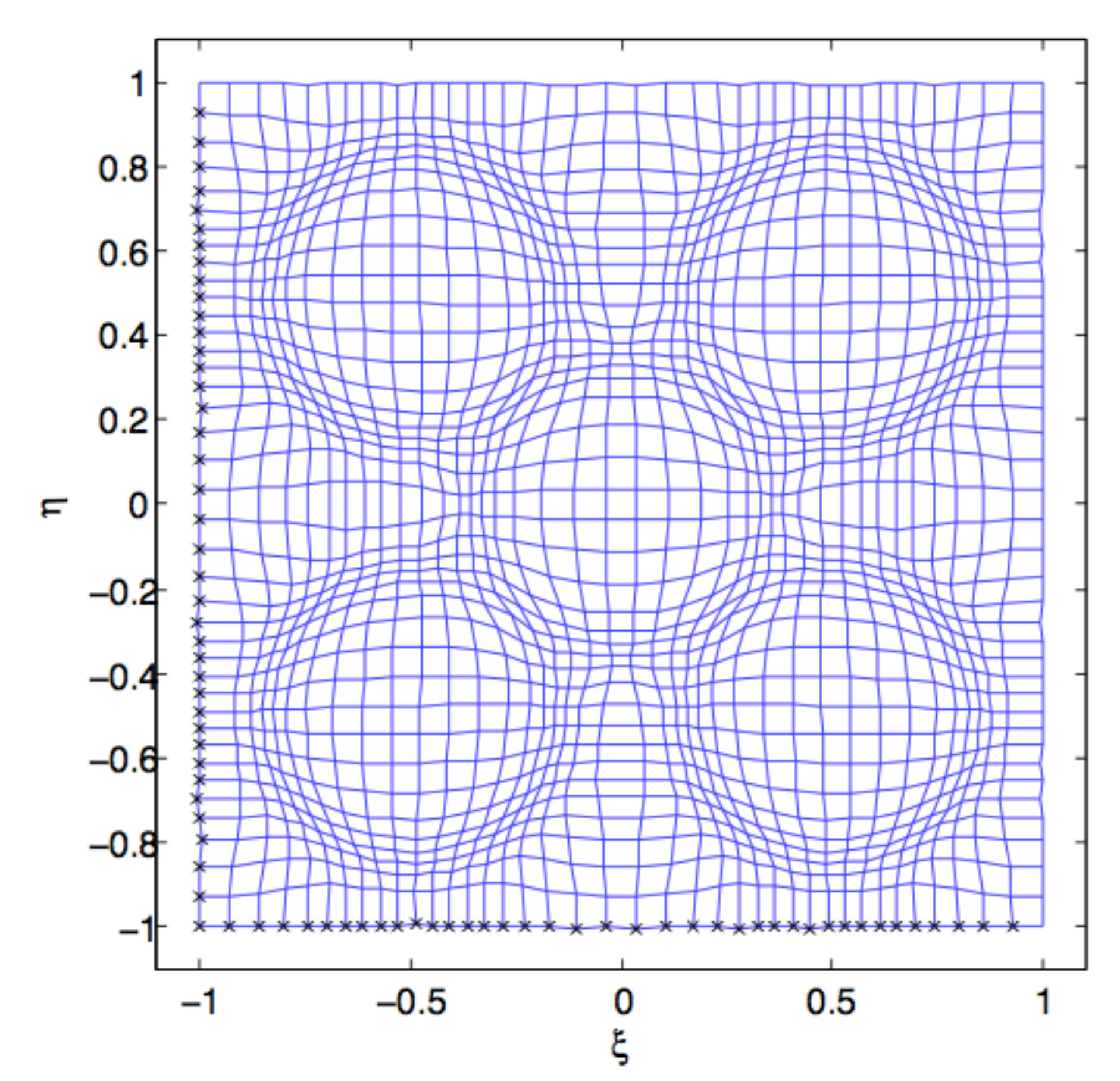}
  \caption{single domain solution}
  \label{fig:ring1}
\end{subfigure}%
\begin{subfigure}{.5\textwidth}
  \centering
  \includegraphics[width=4.5cm,height=4.5cm]{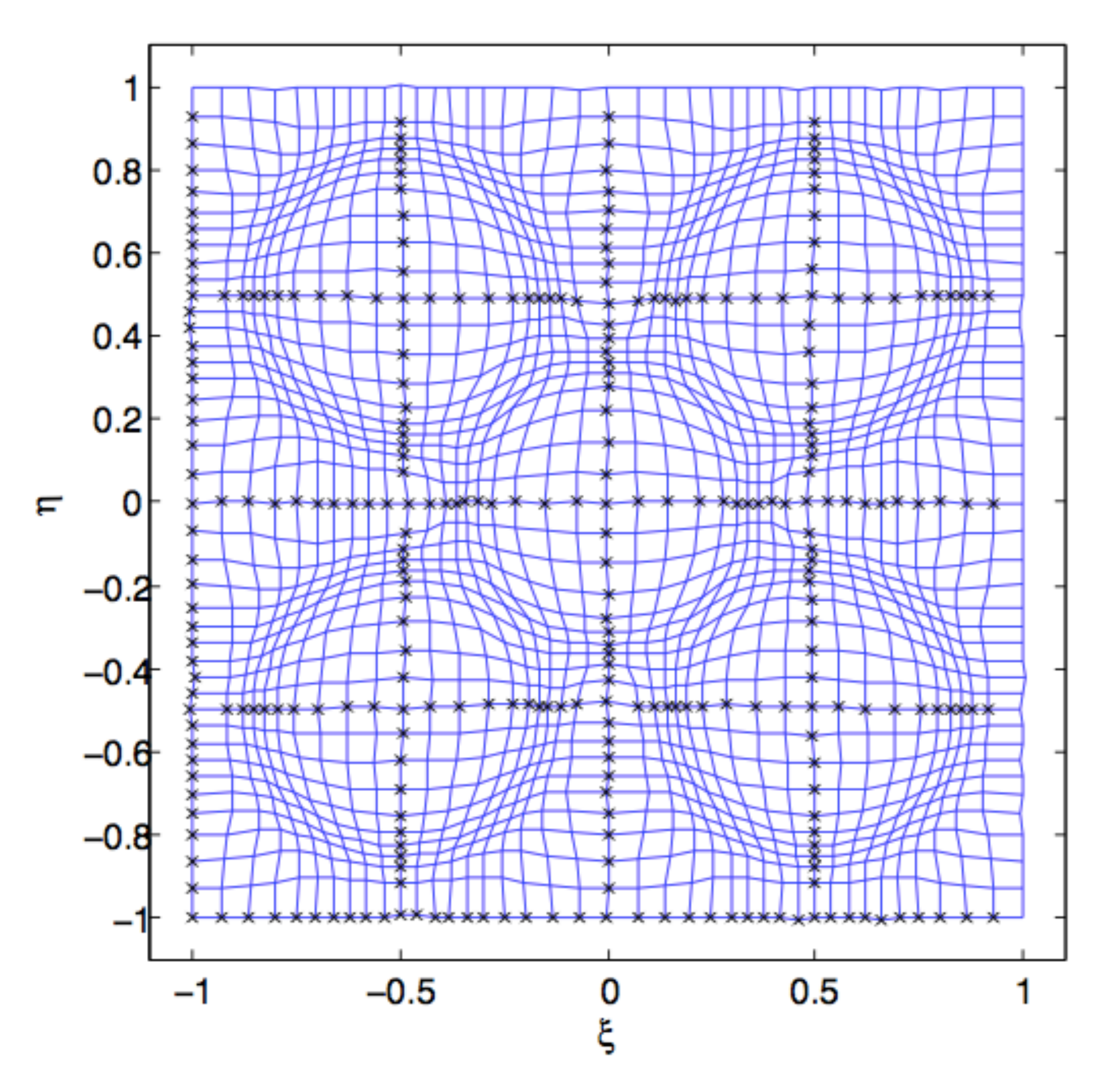}
  \caption{solution with $4\times4$ subdomains }
  \label{fig:ring2}
  \end{subfigure}
    \caption{A comparison of the single domain solution and that obtained for $4\times 4$ subdomains for the five ring problem.}
  \label{fig:ring}
\end{figure}

In order to demonstrate the generation of meshes over periodic domains we first study the case of a prescribed mesh density function. In particular, we re-visit the five ring problem in the form considered in~\cite{haynesbihlo2}. That is, we consider an analytically specified velocity function of the form
\begin{align*}
\begin{split}
 u&=\tanh\left[R\left(x^2+y^2-\frac18\right)\right]+
   \tanh\left[R\left(\left(x-\frac12\right)^2+\left(y-\frac12\right)^2-\frac18\right)\right]\\
   &+\tanh\left[R\left(\left(x-\frac12\right)^2+\left(y+\frac12\right)^2-\frac18\right)\right]+
     \tanh\left[R\left(\left(x+\frac12\right)^2+\left(y-\frac12\right)^2-\frac18\right)\right]\\
   &+\tanh\left[R\left(\left(x+\frac12\right)^2+\left(y+\frac12\right)^2-\frac18\right)\right],
\end{split}
\end{align*}
on the domain $\Omega_{\rm p}=[-1,1[\times[-1,1[$, with $R=30$. A simple arc-length mesh density function
\[
 \rho=\sqrt{1+\alpha(u_x^2+u_y^2)},
\]
with $\alpha=0.2$ is used. The physical domain is discretized using $41\times 41$ grid points, the final integration time is $t=0.05$ with a time step $\Delta t=5\times10^{-4}$. Since this time step is quite small, we could use $\Delta t_s=\Delta t$ for the solution of the discretized SDE~\eqref{eq:WinslowGeneratorSDEDiscretization} as well. Again, $N=10 000$ Monte--Carlo simulations were used. In Fig. {\ref{fig:ring}} the mesh generated on a single domain is compared to the mesh generated by splitting the physical domain into $4\times4$ subdomains.

On inspection these meshes appear to be of similar quality and the computed measures of mesh quality confirm this, with $R_{\max}=1$ and $R_{\rm mean}=1$ for
the $2\times 2$, $3\times 3$ and $4\times 4$ configurations.  The mesh quality is identical to the single domain reference case.



\subsection{Burgers' equation on a periodic domain}

We repeat here the integration of Burgers' equation~\eqref{eq:BurgersEquation2D} but now using a doubly periodic domain of size $\Omega_{\rm p}=[0,2\pi[\times [0,2\pi[$. The initial condition is
\[
 u=u_0+A\sin(x+y-2\pi),
\]
where $u_0=0.75$ and $A=0.5$. The diffusion coefficient is $\nu=0.001$. Again, $41\times41$ grid points are used and the time step of the integration was $\Delta t=0.005$, which coincides with the time step used in the SDE~\eqref{eq:WinslowGeneratorSDEDiscretization}. Only $N=1 000$ Monte--Carlo simulations were used. The results of the integration at times $t=0.85$ and $t=1$ are collected in Table~\ref{tab:BurgersPeriodic}.

\begin{table}[!ht]
\centering
\caption{Mesh qualities for Burgers' equation with periodic boundary conditions.}
\begin{tabular}{|r||r|r||r|r||r|r|}
  \hline
  $t_f$ & $R^{2\times 2}_{\max}$ & $R^{2\times 2}_{\rm mean}$ & $R^{3\times 3}_{\max}$ & $R^{3\times 3}_{\rm mean}$ & $R^{4\times 4}_{\max}$ & $R^{4\times 4}_{\rm mean}$\\
  \hline\hline
  0.85 & 0.99 & 1 & 1 & 1 & 0.99 & 1 \\
  1 & 0.99 & 1 & 0.99 & 1 & 0.99 & 1 \\
  \hline
\end{tabular}
\label{tab:BurgersPeriodic}
\end{table}

The results in Table~\ref{tab:BurgersPeriodic} confirm that the domain decomposition solution leads to meshes with almost the same geometric mesh quality as the single domain reference solution.

\subsection{Shallow Water Equations on a Periodic Domain}

In this final example we solve the system of shallow-water equations in nondimensional form
\begin{align}\label{eq:ShallowWaterEquations}
\begin{split}
& u_t+uu_x+vu_y + h_x = 0,\\
& v_t+uv_x+vv_y + h_y = 0,\\
& h_t+uh_x+vh_y + h(u_x+v_y) = 0,\\
\end{split}
\end{align}
where $(u,v)$ is the two-dimensional velocity field and $h$ is the height of a water column over a constant reference level. For this problem, periodic boundary conditions are considered.

We discretize system~\eqref{eq:ShallowWaterEquations} in computational coordinates using centered
differences for the spatial derivatives and a trapezoidal rule for the time integration. The physical domain is of size $\Omega_{\rm p}=[0,2\pi[ \times [0,2\pi[$, which is discretized using $41\times41$ grid points. The time step of the integration was $\Delta t=0.005$, which again coincides with the time step used for the solution of the SDE~\eqref{eq:WinslowGeneratorSDEDiscretization}. We found that using $N=1 000$ Monte--Carlo simulations gives sufficiently good results.

The initial condition is a pile of water of height $h=12.5$ in the center of the domain over a base level at height $h=10$. This initial condition simulates the breaking of a dam, i.e.\ the evolution of the water level once the walls holding the pile of water have been removed.

\begin{figure}[!ht]
\begin{subfigure}{.5\textwidth}
  \centering
  \includegraphics[width=6cm,height=6cm]{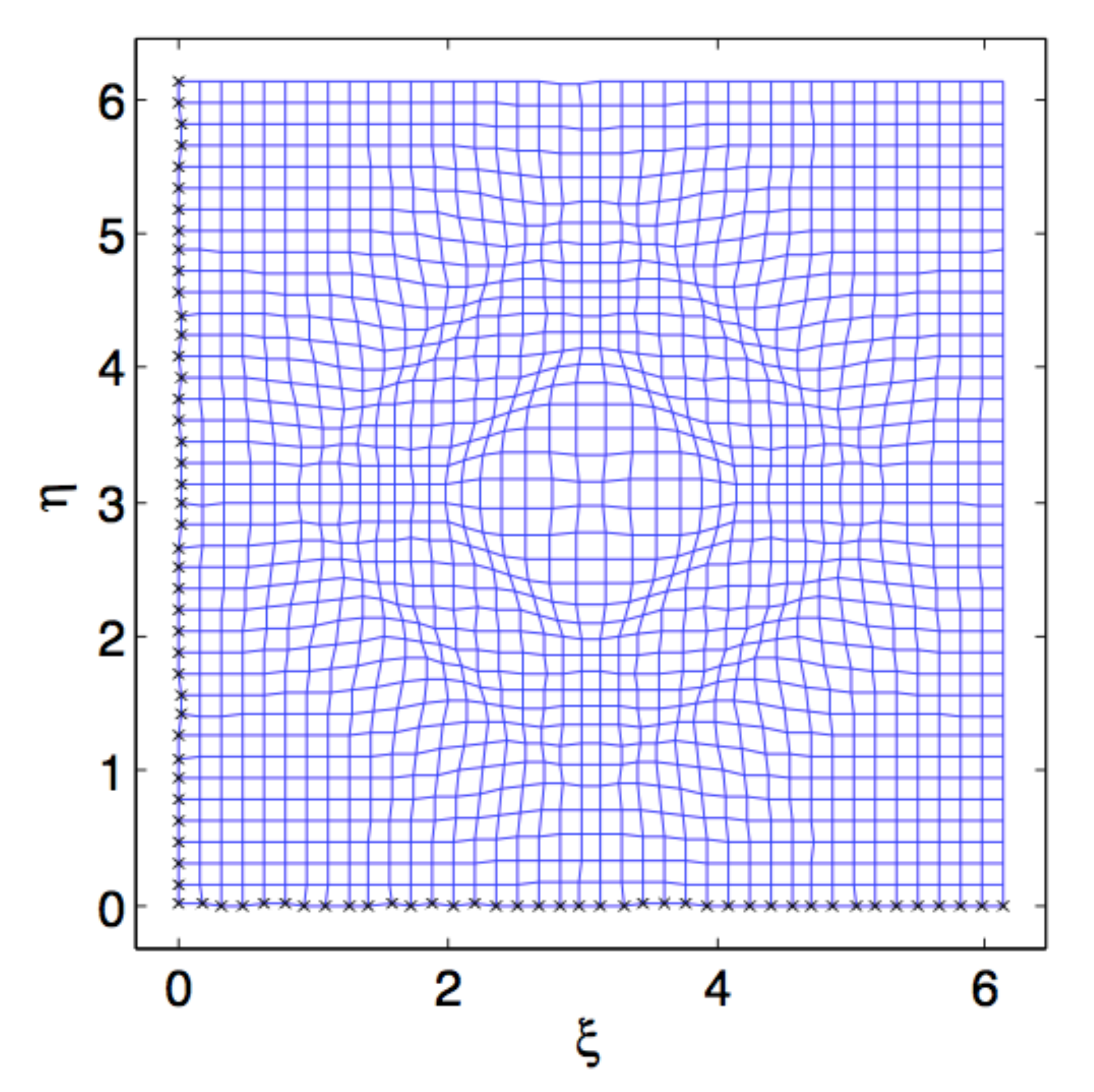}
  \caption{single domain solution}
  \label{fig:sw1}
\end{subfigure}%
\begin{subfigure}{.5\textwidth}
  \centering
  \includegraphics[width=6cm,height=6cm]{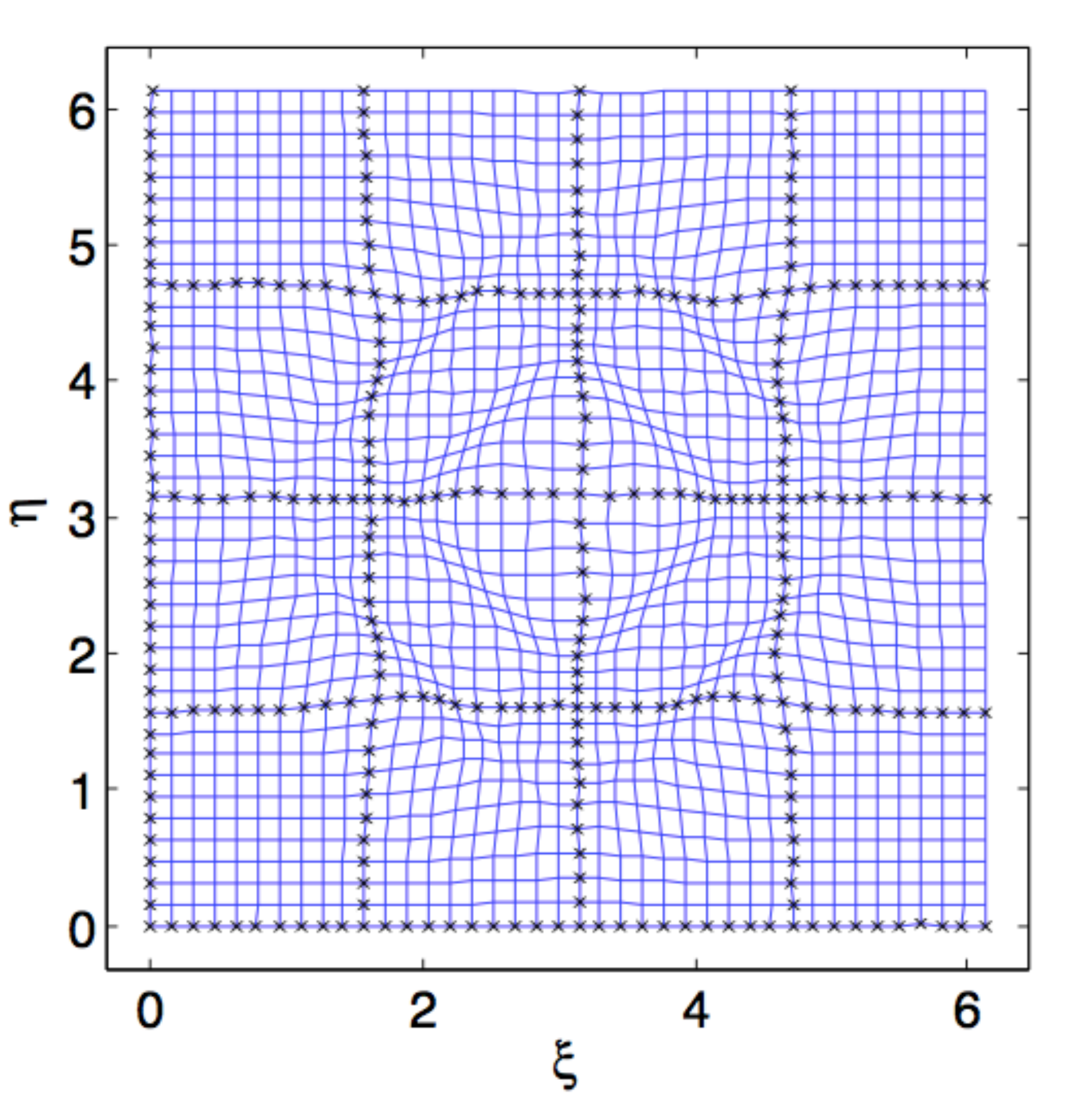}
  \caption{solution with $4\times4$ subdomains }
  \label{fig:sw2}
  \end{subfigure}
    \caption{A comparison of the single domain solution and that obtained for $4\times 4$ subdomains for the shallow water equations at $t=0.15$.}
  \label{fig:sw}
\end{figure}

Since the initial pile of water decays rapidly into gravity waves, we kept the integration times short and report our mesh quality results only at $t=0.05$ and $t=0.15$. A comparison of the mesh for the single domain solution and that computed for $4\times 4$ subdomains can be seen in Fig.~{\ref{fig:sw}}.

As we saw in the previous examples the meshes appear to be of similar quality and this is substantiated by the computed mesh quality measures which can be found in Table~\ref{tab:ShallowWaterPeriodic}. The mesh qualities obtained from the domain decomposition solution are very close to that of the single domain reference solution.

\begin{table}[!ht]
\centering
\caption{Mesh quality for the shallow-water equations with periodic boundary conditions.}
\begin{tabular}{|r||r|r||r|r||r|r|}
  \hline
  $t_f$ & $R^{2\times2}_{\max}$ & $R^{2\times2}_{\mathrm{mean}}$ & $R^{3\times3}_{\max}$ & $R^{3\times3}_{\mathrm{mean}}$ & $R^{4\times4}_{\max}$ & $R^{4\times4}_{\mathrm{mean}}$ \\
  \hline\hline
  $0.05$ & 0.99 & 1 & 0.97 & 1 & 0.99 & 1 \\
  $0.15$ & 0.99 & 1 & 0.99 & 1 & 0.99 & 1 \\
  \hline
\end{tabular}
\label{tab:ShallowWaterPeriodic}
\end{table}

\section{Conclusion}\label{sec:conclude}

Originally, the SDD method has been proposed in~\cite{aceb05a} for the parallel solution of linear elliptic boundary value problems. In~\cite{haynesbihlo2,bihl15a} we have extended the method to the parallel generation of adaptive meshes for prescribed mesh density functions. In this paper we have for the first time demonstrated the use of stochastic domain decomposition for the generation of adaptive moving meshes for time dependent PDE problems.  Furthermore, we have provided the first stochastic domain decomposition method for the generation of  meshes on periodic domains.

Our algorithm hinges on the alternating MP procedure for PDE based mesh generation in higher dimensions.    A fully parallel algorithm would result if the physical PDE was also
solved in parallel using domain decomposition or other approaches.

Indeed the linearization
provided by the MP procedure has the undesirable affect of decoupling the physical solution from the mesh.
This effect
can be reduced using a
$M^kP$ procedure.  This produces a mesh which more closely satisfies the equidistribution principle.
In this case, we approximate $x^{n+1}$ by a sequence of {\em sub-meshes} $x^{n+1,k}$, where $x^{n+1,k+1}$ is obtained from
$x^{n+1,k}$ by using a step of $\Delta t_n/K$ with a  linearized MMPDE.
In this case $\rho$ is approximated by $\rho^{n+1,k}$ obtained by constructing a piecewise linear interpolant of
$(x^n,\rho^n)$ values onto $x^{n+1,k}$.
A $M^{\nu}P$ algorithm, which updates from
$t_n$ to $t_{n+1}$ using variable time steps can also be used.
The time lag between the mesh and physical solution can be reduced
by iterating between the mesh and physical PDE $l$ times,
resulting in a $(MP)^l$ algorithm.   The algorithm
$(MP)^\infty$ is equivalent to the simultaneous solution (if the iteration
converges).   These solution variants are discussed at length in \cite{huan10a}.
The implementation of these variants and a study of their efficacy within the stochastic domain decomposition framework is a topic of
current investigation.

During our experiments we also detected that
fewer Monte--Carlo simulations are needed to generate quality periodic meshes.  This seems to be due to the absence of an exit time test for the periodic case.
Work to fully understand and utilize this is also underway.

\section*{Acknowledgements}

The authors thank Weizhang Huang for helpful discussions. AB is a recipient of an APART Fellowship of the Austrian Academy of Sciences. This research was supported by the Natural Sciences and Engineering Research Council of Canada (NSERC).

{\footnotesize\itemsep=0ex
\bibliography{bibliographyStochasticDD}

\def\cprime{$'$}
\begin{thebibliography}{10}
\providecommand{\url}[1]{\texttt{#1}}
\providecommand{\urlprefix}{URL }
\expandafter\ifx\csname urlstyle\endcsname\relax
  \providecommand{\doi}[1]{doi:\discretionary{}{}{}#1}\else
  \providecommand{\doi}{doi:\discretionary{}{}{}\begingroup
  \urlstyle{rm}\Url}\fi
\providecommand{\eprint}[2][]{\url{#2}}

\bibitem{aceb05a}
Acebr{\'o}n J.A., Busico M.P., Lanucara P. and Spigler R., Domain decomposition
  solution of elliptic boundary-value problems via monte carlo and quasi-monte
  carlo methods, \emph{SIAM J. Sci. Comput.} \textbf{27} (2005), 440--457.

\bibitem{aceb07a}
Acebr{\'o}n J.A. and Spigler R., A new probabilistic approach to the domain
  decomposition method, in \emph{Domain Decomposition Methods in Science and
  Engineering XVI}, Springer, pp. 473--480, 2007.

\bibitem{baines2009}
Baines M.J., Hubbard M.E., Jimack P.K. and Mahmood R., A moving-mesh finite
  element method and its application to the numerical solution of phase-change
  problems, \emph{Commun. Comput. Phys.} \textbf{6} (2009), 595--624.

\bibitem{beckett2001}
Beckett G., Mackenzie J.A. and Robertson M.L., A moving mesh finite element
  method for the solution of two-dimensional stefan problems, \emph{J. Comput.
  Phys.} \textbf{168} (2001), 500--518.

\bibitem{haynesbihlo2}
Bihlo A. and Haynes R.D., Parallel stochastic methods for pde based grid
  generation, \emph{Comput. Math. Appl.} \textbf{68} (2014), 804--820.

\bibitem{bihl15a}
Bihlo A. and Haynes R.D., A stochastic domain decomposition method for time
  dependent mesh generation, in \emph{Domain Decomposition Methods in Science
  and Engineering XXII}, Springer, 2015. In Press.

\bibitem{budd2012}
Budd C.J., Cullen M.J.P. and Walsh E.J., {Monge–Amp\`{e}re based moving mesh
  methods for numerical weather prediction, with applications to the Eady
  problem}, \emph{J.\ Comput. Phys.} \textbf{236} (2013), 247--270.

\bibitem{burchard}
Burchard H.G., Splines (with optimal knots) are better, \emph{Appl. Anal.}
  \textbf{3} (1974), 309--319.

\bibitem{cao99c}
Cao W., Huang W. and Russell R.D., An r-adaptive finite element method based
  upon moving mesh pdes, \emph{J. Comput. Phys.} \textbf{149} (1999), 221--244.

\bibitem{Chen:2012uo}
Chen R. and Cai X.C., {Parallel One-Shot Lagrange--Newton--Krylov--Schwarz
  Algorithms for Shape Optimization of Steady Incompressible Flows}, \emph{SIAM
  J. Sci. Comput.} \textbf{34} (2012), B584--B605.

\bibitem{deboor}
de~Boor C., Good approximation by splines with variable knots, in \emph{Spline
  functions and approximation theory}, Springer, pp. 57--72, 1973.

\bibitem{Haynes:2012}
Gander M.J. and Haynes R.D., Domain decomposition approaches for mesh
  generation via the equidistribution principle, \emph{SIAM J. Numer. Anal.}
  \textbf{50} (2012), 2111--2135.

\bibitem{hayn12a}
Haynes R.D. and Howse A.J.M., Generating equidistributed meshes in 2d via
  domain decomposition, in \emph{Domain Decomposition Methods in Science and
  Engineering XXI}, Springer, pp. 776--797, 2012.

\bibitem{he2012}
He P. and Tang H., An adaptive moving mesh method for two-dimensional
  relativistic magnetohydrodynamics, \emph{Comput. \& Fluids} \textbf{60}
  (2012), 1--20.

\bibitem{huan10a}
Huang W. and Russell R.D., \emph{Adaptive Moving Mesh Methods}, Springer, New
  York, 2010.

\bibitem{Huang2005}
Huang W. and Zhan X., Adaptive moving mesh modeling for two dimensional
  groundwater flow and transport, Amer. Math. Soc., Providence, RI, pp.
  239--252, 2005.

\bibitem{jans03a}
Jansons K.M. and Lythe G.D., Exponential timestepping with boundary test for
  stochastic differential equations, \emph{SIAM J. Sci. Comput.} \textbf{24}
  (2003), 1809--1822.

\bibitem{jin2010}
Jin C., Xu K. and Chen S., A three dimensional gas-kinetic scheme with moving
  mesh for low-speed viscous flow computations, \emph{Adv. Appl. Math. Mech}
  \textbf{2} (2010), 746--762.

\bibitem{kaku44b}
Kakutani S., On {B}rownian motions in $n$-space, \emph{Proc. Imp. Acad. Tokyo}
  \textbf{20} (1944), 648--652.

\bibitem{kaku44a}
Kakutani S., Two-dimensional {B}rownian motion and harmonic functions,
  \emph{Proc. Imp. Acad. Tokyo} \textbf{20} (1944), 706--714.

\bibitem{kara91a}
Karatzas I. and Shreve S.E., \emph{Brownian motion and stochastic calculus},
  vol. 113 of \emph{Graduate Texts in Mathematics}, Springer, New York, 1991.

\bibitem{lang}
Lang J., Cao W., Huang W. and Russell R.D., A two-dimensional moving finite
  element method with local refinement based on a posteriori error estimates,
  \emph{Appl. Numer. Math.} \textbf{46} (2003), 75--94.

\bibitem{mils04a}
Milstein G. and Tretyakov M., \emph{Stochastic numerics for mathematical
  physics}, Springer, Berlin Heidelberg, 2004.

\bibitem{okse10a}
{\O}ksendal B., \emph{{Stochastic Differential Equations: An Introduction with
  Applications}}, Springer, Heidelberg, 2010.

\bibitem{pres07a}
Press W.H., Teukolsky S.A., Vetterling W.T. and Flannery B.P., \emph{Numerical
  recipes 3rd edition: {T}he art of scientific computing}, Cambridge University
  Press, Cambridge, UK, 2007.

\bibitem{quan2011}
Quan S., Simulations of multiphase flows with multiple length scales using
  moving mesh interface tracking with adaptive meshing, \emph{J. Comput. Phys.}
  \textbf{230} (2011), 5430--5448.

\bibitem{tan08a}
Tan Z., Lim K. and Khoo B., An adaptive moving mesh method for two-dimensional
  incompressible viscous flows, \emph{Commun. Comput. Phys.} \textbf{3} (2008),
  679--703.

\bibitem{white}
White Jr A.B., On selection of equidistributing meshes for two-point
  boundary-value problems, \emph{SIAM J. Numer. Anal.} \textbf{16} (1979),
  472--502.

\bibitem{wins66a}
Winslow A.M., Numerical solution of the quasilinear {P}oisson equation in a
  nonuniform triangle mesh, \emph{J. Comput. Phys.} \textbf{1} (1966),
  149--172.

\bibitem{yuan2012}
Yuan Y., The characteristic finite element alternating-direction method with
  moving meshes for the transient behavior of a semiconductor device,
  \emph{Int. J. Numer. Anal. Model.} \textbf{9} (2012), 86--104.

\bibitem{zhang2011}
Zhang Y. and Tang T., Simulating three-dimensional free surface viscoelastic
  flows using moving finite difference schemes, \emph{Numer. Math. Theory
  Methods Appl.} \textbf{4} (2011), 92--112.

\bibitem{zhang}
Zhang Z. and Tang T., \emph{An adaptive mesh redistribution algorithm for
  convection-dominated problems}, Department of Mathematics, Hong Kong Baptist
  University, 2002.

\end{thebibliography}
}

\end{document}